\documentclass[english,a4paper,12pt]{article}
\usepackage{amsmath}
\usepackage{amsfonts}
\usepackage{amssymb}
\usepackage{babel}
\usepackage{latexsym}
\newtheorem{theorem}{Theorem}[section]

\newtheorem{proposition}[theorem]{Proposition}
\newtheorem{definition}[theorem]{Definition}
\newtheorem{definitions}[theorem]{Definitions}
\newtheorem{remark}[theorem]{Remark}
\newtheorem{remarks}[theorem]{Remarks}
\newtheorem{example}[theorem]{Example}
\newtheorem{examples}[theorem]{Examples}

\newcommand{\erre}{\mbox{$\mathbb{R}$}}
\newcommand{\enne}{\mbox{$\mathbb{N}$}}

\title{\Large \bf Korovkin-type theorems for abstract modular
convergence}
\author{A. Boccuto \thanks{
Dipartimento di Matematica e Informatica, University of Perugia,
via Vanvitelli 1, I-06123 Perugia, Italy, Email: antonio.boccuto@unipg.it, boccuto@yahoo.it
(Corresponding author)}
\and  X. Dimitriou\thanks{ Department of Mathematics, University of Athens,
Panepistimiopolis, Athens 15784, Greece,
Email xenofon11@gmail.com
 \newline
 \textit{2010 A. M. S. Subject Classifications:}
40A35, 41A35, 46E30. \newline {\it Key words}:
axiomatic convergence,
modular space, (quantitative)
Korovkin theorem, rate of convergence, filter convergence,
(triangular) statistical convergence. }}
\date{}
\begin{document}
\maketitle
\vskip0,5cm
\begin{abstract}
We give some Korovkin-type theorems on 
convergence and estimates of rates of approximations
of nets of functions, satisfying suitable axioms, whose particular cases
are filter/ideal convergence, almost convergence and triangular 
$A$-statistical convergence, where $A$ is a non-negative summability
method. Furthermore, we give some applications to Mellin-type
convolution and bivariate Kantorovich-type discrete operators.
\end{abstract}
\section{Introduction}
The classical Bohman-Korovkin theorem is a result which
yields uniform convergence in the space ${\mathcal C}([a,b])$ of all continuous real-valued functions defined on the compact subinterval $[a,b]$ of the real line, for a net $(T_w)_w$ of positive linear operators   on $C([a,b])$, with the only hypothesis  of convergence on the test functions $1$,  $x$, $x^2$ (see also \cite{bohman, devore,
korovkin1, korovkin}).  
There have been several extensions
the Korovkin theorem to the context of abstract functional spaces.
For a related literature, see for instance 
\cite{donner, kittowulbert, renaud, wulbert} for 
the case of $L^p$-spaces,
\cite{maligranda, soardi} for Orlicz spaces and 
\cite{bbdmkorovkin, bm2009, belen} for general modular spaces.
There have been also several studies about Korovkin-type 
theorems in the setting of
convergence generated by summability matrices, statistical and 
filter convergence (see for example \cite{AD2, 
duman, dumanorhan, Kdd2, 
demirciorhan2}) and with respect to 
``triangular $A$-statistical convergence'', namely an extension 
of statistical convergence for double sequences of positive linear 
operators, where $A$ is a suitable non-negative regular matrix (see also \cite{bbdmo2, dobbm}).

In this paper we deal with Korovkin-type results about 
convergence and quantitative theorems, and
estimates of rates of approximation 
with respect to abstract convergences satisfying suitable axioms
(see also \cite{bbdmkorovkin, BC,
KURATOWSKI}), including as particular cases
convergence generated by summability (double infinite) matrices,
filter convergence and almost convergence, which is not generated
by any filter (see also \cite{bdpPL}), and we 
consider the general case of a net of operators, acting on an abstract 
modular function space 
generated by a modular, extending earlier results proved in
\cite{bbdmo2, bbdmkorovkin, bdrates, dobbm, duman}. 
Note that, in the literature, it is often dealt with 
Korovkin-type theorems with respect to some of the above
mentioned convergences, but without a general approach 
containing all of them. Our general results unify various previous
theorems.

Furthermore, we show that the rates investigated in 
\cite{bbdmo2, bdrates, dobbm} 
are particular cases of those treated here with respect to the
``axiomatic convergence'', we present
some examples of convergences and give
comparison results between triangular convergence
generated by a summability matrix method and filter
convergence. In particular we see that, in general, 
the $\Psi$-$A$-statistical convergence studied in \cite{bbdmo2}
and the filter convergence are such that neither contains the other.
We give also some applications to moment kernels and bivariate 
Kantorovich-type discrete operators 
(for recent studies and developments, see also \cite{AL, 
AV, bbdm, 
barman2011, barman2012, barman2013bis, 
BMV, belen, bcs, BD1, BD2, bdrates, bdbook, V}).

\section{Preliminaries}
Let $(W, \succeq)$ be a directed set, and let us consider an axiomatic
abstract convergence on $W$, defined as follows
(see also \cite{bbdmkorovkin, BC, KURATOWSKI}).
\begin{definition}\label{axiomaticconvergence} \rm
Let ${\mathcal T}$ be the set of all
real-valued nets $(x_w)_{w\in W}$.
A \textit{convergence} is a pair $({\mathcal S},{ \ell})$, where
${\mathcal S}$ is a linear subspace of ${\mathcal T}$ and
${\ell}:{\mathcal S} \to \erre$ is a function, satisfying the following
axioms:
\begin{description}
\item[{\rm (a)}] ${\ell}((a_1 \, x_w + a_2 \, y_w)_w)=
a_1 \, {\ell} ( (x_w)_w)+ a_2 \, {\ell} ( (y_w)_w)$
for every pair of nets $(x_w)_w$, $(y_w)_w \in {\mathcal S}$ and
for each $a_1$, $a_2 \in \erre$ (linearity).
\item[{\rm (b)}] If $(x_w)_w$, $(y_w)_w
\in {\mathcal S}$ and there is $w^*\in W$ with
$x_w \leq y_w$ for every $w\succeq w^*$,
then ${\ell}((x_w)_w) \leq {\ell}((y_w)_w)$
(monotonicity).
\item[{\rm (c)}] If $(x_w)_w$ is such that there is $w_*\in W$
with $x_w=l$ whenever $w\succeq w_*$, then $(x_w)_w \in {\mathcal
S}$ and ${\ell}((x_w)_w)=l$.
\item[{\rm (d)}] If $(x_w)_w \in {\mathcal S}$, then $(|x_w|)_w \in {\mathcal S}$
and ${\ell}((|x_w|)_w)=|{\ell}((x_w)_w)|$.
\item[{\rm (e)}]
Let $(x_w)_w$, $(y_w)_w$, $(z_w)_w,$ satisfying
$(x_w)_w$, $(z_w)_w \in {\mathcal S}$, ${\ell}((x_w)_w)={\ell}
((z_w)_w)$ and suppose that there is $\overline{w}\in W$
with $x_w \leq y_w \leq z_w$ for every $w\geq \overline{w}$. Then $(y_w)_w \in {\mathcal S}$.
\end{description}
\end{definition}
Note that ${\mathcal S}$ is the space of all
convergent nets, $\ell$ will be the ``limit'' according to this
approach, and we will denote by the symbol $\displaystyle{
(\ell)\lim_w x_w}$ the quantity ${\ell}((x_w)_w)$.


We now give the axiomatic definition of the operators ``limit superior''
and ``limit inferior'' related with
a convergence $({\mathcal S},{\ell})$ 
(see also \cite{bbdmkorovkin}).
\begin{definition}\label{axiomaticlimsupliminf} \rm
Let ${\mathcal T}$, ${\mathcal S}$ be as above.
We define two functions
$\overline{\ell}$,
$\underline{\ell}:{\mathcal T} \to \widetilde{\erre}$, satisfying the
following axioms:
\begin{description}
\item[{\rm (f)}] If $(x_w)_w$, $(y_w)_w
\in {\mathcal T}$,
then $\underline{\ell} ( (x_w)_w) \leq \overline{\ell} ( (x_w)_w)$
and $\overline{\ell} ( (x_w)_w) = - \underline{\ell} ( (-x_w)_w)$.
\item[{\rm (g)}] If $(x_w)_w \in {\mathcal T}$, then

(i) $\overline{\ell}((x_w + y_w)_w) \leq
\overline{\ell} ( (x_w)_w)+  \overline{\ell} ( (y_w)_w)$
(subadditivity);

(ii) $\underline{\ell}((x_w + y_w)_w) \geq
\underline{\ell} ( (x_w)_w)+  \underline{\ell} ( (y_w)_w)$
(superadditivity).
\item[{\rm (h)}] If $(x_w)_w$, $(y_w)_w
\in {\mathcal T}$ and $x_w \leq y_w$ definitely,
then $\overline{\ell}((x_w)_w) \leq \overline{\ell}((y_w)_w)$ and
$\underline{\ell}((x_w)_w) \leq \underline{\ell}((y_w)_w)$
(monotonicity).
\item[{\rm (j)}] A net $(x_w)_w \in {\mathcal T}$
belongs to ${\mathcal S}$  if and only if
$\overline{\ell}((x_w)_w) = \underline{\ell}((x_w)_w)$.
\end{description}
We will denote by the symbols $\displaystyle{
(\ell)\limsup_w x_w}$ and $\displaystyle{
(\ell)\liminf_w x_w}$
the quantities $\overline{\ell}((x_w)_w)$ and
$\overline{\ell}((x_w)_w)$, respectively.
\end{definition}
We define two tools, in order to
``compare'' two nets belonging to ${\mathcal S}$.
\begin{definition} \rm  Let $(x_w)_w$, $(y_w)_w
\in {\mathcal S}$ with $\displaystyle{(\ell)\lim_w x_w=
(\ell)\lim_w y_w=0}$
and $y_w\neq 0$ for every $w\in W$. We say that $x_w=o(y_w)$ iff
$\displaystyle{(\ell)\lim_w \frac{|x_w|}{|y_w|}=0}$, and that
$x_w=O(y_w)$ iff
$\displaystyle{(\ell)\limsup_w \frac{|x_w|}{|y_w|}\in\mathbb{R}}$.
\end{definition} \vspace{3mm}


Let $G=(G,d)$ be a metric space,
${\mathcal B}$
be the $\sigma$-algebra of all Borel subsets of $G$, and $\mu$ be a
positive finite regular measure defined on ${\mathcal B}$.
Let $L^0(G)$ be the space of all
real-valued $\mu$-measurable functions on $G$ with identification
up to sets of measure $\mu$ zero, ${\mathcal C}_b(G)$ be the space of all real-valued
continuous and bounded functions on $G$, ${\mathcal C}_c(G)$ be
the subspace of ${\mathcal C}_b(G)$ of all
functions with compact support on $G$ and Lip$(G)$ be
the space of all real-valued
Lipschitz functions on $G$.

We now recall the notion of modular space
(see also \cite{BMV}).
\begin{definitions}\label{modulardef} \rm
(a) A functional $\rho:L^0(G) \to \widetilde{\mathbb{R}^+_0}$ is called
a \textit{modular} on $L^0(G)$ iff it satisfies the following conditions:

i) $\rho[f]=0 \Longleftrightarrow f=0$ $\mu$-almost everywhere on $G$;

ii) $\rho[-f]=\rho[f]$ for every $f \in L^0(G)$;

iii) $\rho[a f + b g] \leq \rho[f] + \rho[g]$ whenever
$f$, $g \in L^0(G)$ and $a \geq 0$, $b \geq 0$ with
$a + b =1$.

(b) A modular $\rho$ is said to be \textit{convex} iff it satisfies conditions
i), ii) and

iii') $\rho[a f + b g] \leq a \rho[f] + b \rho[g]$ for all
$f$, $g \in L^0(G)$ and for every $a$, $b \geq 0$ with
$a + b=1.$

(c)
Let $Q \geq 1$ be a real constant.
We say that a modular $\rho$ is \textit{$Q$-quasi semiconvex} if
$\rho[a \, f] \leq Q \, a \, \rho[Q \, f]$ for all $f \in L^0(G)$, $f \geq 0$ and
$0 < a \leq 1$ (see also \cite{bm2009}).


(d) A modular $\rho$ is \textit{monotone} iff $\rho[f]
\leq \rho[g]$ for all $f$, $g \in L^0(G)$ with
$|f|\leq |g|$.

(e) A modular $\rho$ is \textit{finite} iff $\chi_A$ (the characteristic
function associated with $A$) belongs to $L^{\rho}(G)$ whenever
$A \in {\mathcal B}$ with $\mu(A) < + \infty$.

(f) A modular $\rho$ is \textit{strongly finite} iff $\chi_A$
belongs to $E^{\rho}(G)$ for each $A \in {\mathcal B}$
with $\mu(A) < + \infty$.

(g) A modular $\rho$ is said to be \textit{absolutely continuous} iff
there is a positive constant $a$ with the property: for all
$f \in L^0(G)$ with $\rho[f]< + \infty$,

i) for each $\varepsilon > 0$ there exists a set $A \in {\mathcal B}$
with $\mu(A) < + \infty$ and $\rho[a f \chi_{G \setminus A}] \leq \varepsilon $,

ii) for every $\varepsilon >0$ there is a $\delta >0$
with $\rho[a f \, \chi_{B}] \leq \varepsilon$ for every
$B \in {\mathcal B}$ with $\mu(B) < \delta$.

(h) The \textit{modular space} $L^{\rho}(G)$
generated by $\rho$ is $$L^{\rho}(G)=\{f \in L^0(G):
\lim_{\lambda \to 0^+ } \rho[\lambda f]=0 \},$$ 
where the limit is intended in the usual sense, and
the \textit{space of the finite elements of} $L^{\rho}(G)$ is
$$E^{\rho}(G)=\{f \in L^{\rho}(G):
 \rho[\lambda f] < + \infty \text{ for all } \lambda >0 \}.$$
\end{definitions}
\begin{example}\label{orlicz} \rm Let
$\Phi$ be the set of all continuous non-decreasing functions
$\varphi: \mathbb{R}^+_0 \to \mathbb{R}^+_0$ with $\varphi(0)=0$,
$\varphi(u) >0$ for any $u >0$ and $\displaystyle{\lim_{u \to + \infty}
\varphi(u)=+\infty}$ in the usual sense,
and let $\widetilde{\Phi}$ be the set of all convex functions belonging
to $\Phi$.

For every $\varphi \in \Phi$ (resp. $\widetilde{\Phi})$,
the functional $\rho^{\varphi}$ defined by
\begin{eqnarray}\label{modularorlicz}
\rho^{\varphi}[f]=\int_G \varphi(|f(s)|) \, d\mu(s),
\quad f \in L^0(G),
\end{eqnarray}
is a modular (resp. convex modular) on $L^0(G)$ and
$$L^{\varphi}(G):=\{ f \in L^0(G): \rho^{\varphi}[\lambda f]< + \infty
\text{ for some } \lambda >0 \}$$ is the
\textit{Orlicz space} generated by $\varphi$,
and satisfies all above properties (see also \cite{BMV}).
\end{example}
We now define the modular and strong convergences in the context of
the axiomatic convergence in Definition \ref{axiomaticconvergence}
(for the classical case and filter convergence see \cite{BMV}
and \cite{bbdmkorovkin, bbdm, bcs} respectively).

A net $(f_w)_w$ of functions in $L^{\rho}(G)$
is \textit{$(\ell)$-modularly convergent} to
$f \in L^{\rho}(G)$ if there is a $\lambda >0$ with
$$(\ell)\lim_w \rho[\lambda(f_w-f)]=0 .$$
A net $(f_w)_w$ in $L^{\rho}(G)$ is
\textit{$(\ell)$-strongly convergent} to
$f \in L^{\rho}(G)$  if
$$(\ell)\lim_w \rho[\lambda(f_w-f)]=0  \text{  for
every   } \lambda  >0.$$

Given a subset ${\mathcal A} \subset L^{\rho}(G)$ and $f \in L^{\rho}(G)$, we say that
$f \in \overline{\mathcal {A}}$ (that is, $f$ is {in the modular closure of ${\mathcal A}$})
if there is a sequence $(f_k)_k$ 
in ${\mathcal A}$, modularly convergent to $f$ with respect
to the ordinary convergence.

We recall the following
\begin{proposition}\label{ily}
{\rm (see also \cite[Theorem 1]{ilaria})}
Let $\rho$ be a monotone, strongly finite
and absolutely continuous modular on $L^0(G)$.
Then $\overline{{\mathcal C}_c(G)}=L^{\rho}(G)$ with respect to the
modular convergence in the ordinary sense.
\end{proposition}
\section{The Korovkin theorem}
We consider some kinds of rates of
approximation associated with the Korovkin theorem in the context of
modular convergence. For
technical reasons, we sometimes suppose that
$(G,d)$ satisfies the following property:
\begin{description}
\item[{\rm H*)}] For every $n\in\mathbb{N}$ and $s$, $t\in G$, with $s \neq t$, there are
$n+1$ points $x_i$, $i=0, \ldots, n+1$, such that $s=x_0$, 
$t=x_{n+1}$ and $\displaystyle{
d(x_i,x_{i+1})\leq \frac1n \, d(s,t)}$ for each $i=0, \ldots, n$.
\end{description} Some examples of spaces satisfying
condition H*) are the Euclidean multidimensional space
$\mathbb{R}^N$ 
endowed with the usual metric and the space
${\mathbb{R}}^{\Lambda}$ equipped with the sup-norm, where
$\Lambda$ is any abstract nonempty set (see also \cite{bbdmo2}).

For every $f \in {\mathcal C}_b(G)$ and $\delta >0$, let
$$\omega (f; \delta):= \sup \{|f(s)-f(t)|:s,t\in G, d(s,t) \leq \delta\}$$
be the usual {modulus of continuity} of $f$.
Note that $\omega(f;\delta)$ is an increasing function of $\delta$,
$|f(s)-f(t)|\leq\omega(f;d(s,t))$ for each $s$, $t\in G$,
$\omega(f;\delta) \leq 2 \, M$ for every $\delta$, where
$\displaystyle{M=\sup_{t\in G}|f(t)|}$, and
\begin{eqnarray}\label{modulus}
\omega(f;\gamma \,\delta) \leq (1+\gamma) \, \omega(f; \delta)
\end{eqnarray} for every $\gamma$, $\delta >0$
(see also \cite{bbdmo2}).

Let $T$ be a net of linear
operators $T_w: {\mathcal D} \to L^0(G)$, $w\in W$, with ${\mathcal C}_b(G)
\subset {\mathcal D} \subset L^0(G)$. Here the set ${\mathcal D}$ is the domain of the
operators $T_w$.

We say that the net $T$, together with the modular $\rho$,
satisfies \textit{property $(\rho)$-$(*)$} iff there exist a subset $X_{{T}}
\subset {\mathcal D} \cap L^{\rho}(G)$ with ${\mathcal C}_b(G)
\subset X_{{T}}$
and an $E>0$
with $T_w f \in L^{\rho}(G)$ for any $f \in X_{{T}}$ and $w
\in W$, and $\displaystyle{(\ell)
\limsup_w \rho[\tau (T_w f) ] \leq E \, \rho[\tau f]}$
for every $f \in X_{{T}}$ and $\tau >0$.

Some examples of operators satisfying property $(\rho)$-$(*)$
can be found in \cite{bm2009}.

Let $e_r$ and $a_r$, $r=0, \ldots, m$, be functions in
${\mathcal C}_b(G)$, with
$e_0(t):=1$ for every $t \in G$. Let us define
\begin{eqnarray}\label{pi}
P_s(t):= \sum_{r=0}^m a_r(s) e_r(t), \quad s,t \in G,
\end{eqnarray} and assume that
\begin{description}
\item[{\rm (P1)}] $P_s(s)=0$ for all $s \in G$;
\item[{\rm (P2)}] there is a $C_1>0$
with $P_s(t)\geq C_1 \, d(s,t)$
whenever $s$, $t\in G$.
\end{description}
\begin{examples}\label{esempi} \rm ~

(a) Let $G = I^m$ be endowed with the usual norm $\|\cdot\|_2$, 
where $I\subset \mathbb{R}$ is a connected set, and $\phi:I\to \mathbb{R}$ be
monotone, continuous and such that
$\phi^{-1}$ is Lipschitz on $\phi(I)$.
Examples of such functions are 
$\phi(t)=t$ or $\phi(t)=e^t$, 
where $I=[a,b] \subset \mathbb{R}$.

For every $t=(t_1, \ldots, t_m) \in G$ set $e_i(t):=\phi(t_i)$, $i=1, \ldots, m$,
and $$\displaystyle{e_{m+1}(t):=\sum_{i=1}^m [\phi(t_i)]^2}.$$

For each $s=(s_1, \ldots, s_m) \in G$ put 
$\displaystyle{a_0(s):=\sum_{i=1}^m [\phi(s_i)]^2}$,
$a_i(s)=-2\phi(s_i)$, $i=1, \ldots, m$, and $a_{m+1}(s)\equiv 1$.
We get:
\begin{eqnarray*}
P_s(t)&:=& \sum_{i=0}^{m+1} a_i(s) e_i(t)
=\sum_{i=1}^m [\phi(s_i)-\phi(t_i)]^2.
\end{eqnarray*}
It is not difficult to see that (P1) and (P2) are satisfied.

(b)
Let $G=[a,b]$ with $0 < a < b < \pi/2$, $e_1(t)=\cos t$, $e_2(t)=\sin t$,
$t \in G$. Set $a_0(s)\equiv 1$, $a_1(s)=-\cos s$, $a_2(s)=-\sin s$, $s \in G$. For all
$s$, $t \in G$ we get:
\begin{eqnarray*}
P_s(t)=1-\cos s \cos t - \sin s \sin t=1-\cos(s-t).
\end{eqnarray*}
It is not difficult to check that (P1) and (P2) are fulfilled
(see also \cite{bbdmkorovkin}).
\end{examples}
From now on, we suppose that $e_r \in L^\rho(G),$ $r=0,1,\ldots,
m$. Note that this assumption is fulfilled, for example,
when $G$ is 
a space of finite measure $\mu$.

We now state the following theorems, whose proofs are 
analogous to those of \cite[Theorem 4.2]{bbdmkorovkin} and
\cite[Theorem 4.3]{bbdmkorovkin}, respectively.
\begin{theorem}\label{korovkin1}
Let $\rho$ be a strongly finite, monotone and $Q$-quasi semiconvex modular.
Assume that $e_r$ and $a_r$, $r=0, \ldots, m$, satisfy 
\rm (P1) \em and \rm (P2). \em
Let $T_w$, $w\in W$, be a net of positive
linear operators having property $(\rho)$-$(\ast)$. If
$(T_w e_r)_w$
is $(\ell)$-modularly convergent to $e_r$
in $L^{\rho}(G)$ for each
$r=0, \ldots, m$, then $(T_w f)_w$
is $(\ell)$-modularly convergent
to $f$ in $L^{\rho}(G)$ for every $f \in {\mathcal C}_c(G)$.

If $(T_w e_r)_w$ is $(\ell)$-strongly convergent to $e_r$,
$r=0, \ldots, m$ in $L^{\rho}(G)$, then $(T_w f)_w$
is $(\ell)$-strongly convergent to $f$ in $L^{\rho}(G)$ for every $f
\in {\mathcal C}_c(G).$
\end{theorem}
\begin{theorem}\label{korovkin2}
Let $\rho$ be a monotone, strongly finite,
absolutely continuous and $Q$-quasi semiconvex modular on $L^0(G)$,
and $T_w$, $w\in W$ be a net of positive
linear operators satisfying $(\rho)$-$(*)$. If $(T_w e_r)_w$
is $(\ell)$-strongly convergent to $e_r$,
$r=0, \ldots, m$ in $L^{\rho}(G)$, then $(T_w e_r)_w$
is $(\ell)$-modularly convergent to $f$
in $L^{\rho}(G)$ for every $f \in L^{\rho}(G) \cap {\mathcal D}$
with $f - {\mathcal C}_b(G) \subset X_{{T}}$, where ${\mathcal D}$
and $X_{{T}}$ are as above.
\end{theorem}

Now we present some estimates on rates of approximation
for abstract Korovkin-type theorems.
Let $\Xi$ be the family of all nets $\xi_w$, $w\in W$, with $\xi_w
\neq 0$ for each $w\in W$ and
$\displaystyle{(\ell)\lim_w \xi_w}=0$.

\begin{theorem}\label{ratesdouble}
Let $Q \geq 1$, $\rho$ be a 
monotone, strongly finite and $Q$-quasi semiconvex modular,
$T_w$, $w\in W$, be a net of positive linear operators and $\Xi$ be as above.
For every $w\in W$,
let $\xi^r_w\in \Xi$, $r=0, \ldots, m$,
and set $\xi_w:=\max \{
\xi^r_w$: $r=0, \ldots, m\}$. If $\gamma >0$ is such that
$\rho[\gamma(T_w e_r - e_r)]=o(\xi^r_w)$ for each
$r=0,\ldots, m$, then for every $f\in {\mathcal C}_c(G)\cap \text{Lip}
(G)$ there exists a positive real number $\tau$ 
with $\rho[\tau(T_w f - f)]=o(\xi_w)$.

A similar result holds also when $o$ is replaced by $O$.
\end{theorem}
{\bf Proof:} We now prove only the result concerning $o$, 
since the assertion involving $O$ is analogous.
Choose $f \in {\mathcal C}_c(G)\cap \text{Lip}(G)$,
let $\displaystyle{M:=1+\sup_{t\in G}|f(t)|}$ and
let $C_2$ be a strictly positive
Lipschitz constant, associated with $f$.
We get
%
$$|f(s)-f(t)| \leq C_2 \, d(s,t) \leq
{C_1}^{-1} \, C_2 P_s(t)
\quad \text{  for every  } s,t \in G,$$ namely
\begin{eqnarray}\label{star}
-  {C_1}^{-1} C_2 P_s(t) \leq f(s)-f(t) \leq
{C_1}^{-1} C_2 P_s(t) \quad \text{  for  each  } s,t \in G.
\end{eqnarray}
By applying $T_w$ to (\ref{star}), since $T_w$ is linear and positive
we have
\begin{eqnarray*}
- {C_1}^{-1} C_2 (T_w P_s) (s) \leq
f(s) (T_w e_0) (s)
-(T_w f) (s) \leq
{C_1}^{-1} C_2 (T_w P_s) (s),
\end{eqnarray*} and hence
\begin{eqnarray}\label{prestimab}
|(T_w f)(s) - f(s) | &\leq& |(T_w f)(s) - f(s) (T_w e_0) (s)| +
| f(s) (T_w e_0) (s) -f(s)| \leq \nonumber
\\ &\leq&  {C_1}^{-1} C_2
(T_w P_s) (s) + M |(T_w e_0)(s) -e_0(s)|
\end{eqnarray} for each $s\in G$ and $w\in W$.

Let now $\gamma >0$ be as in the hypotheses and choose a positive
real number $\tau$, with
$$\tau \leq \min \Bigl\{\frac{\gamma \, C_1}{2 \, C_2 (m+1)N\, Q^2},
\frac{\gamma}{2\,M\,Q^2}\Bigr\}.$$ By applying the modular $\rho$,
from (\ref{prestimab}) we get
\begin{eqnarray}\label{modular}
\rho[\tau(T_w f - f)] \leq \rho[ 2 \, \tau
{C_1}^{-1} C_2 (T_w P_{(\cdot)}) (\cdot)] +
\rho [2 \, \tau \, M (T_w e_0 -e_0)]
\end{eqnarray} whenever $w\in W$. By $Q$-quasi semiconvexity of
$\rho$ we have
\begin{eqnarray*}
\rho [2 \, \tau \, M (T_w e_0 -e_0)] &\leq& \rho 
\Bigl[\frac{\gamma}{Q^2}
(T_w e_0 - e_0) \Bigr] \leq Q \, \frac1Q \rho \Bigl[ Q \,
\frac{\gamma}{Q}(T_w e_0 -e_0)\Bigr]=\\&=&\rho [\gamma (T_w e_0 
-e_0)].
\end{eqnarray*} Hence,
\begin{eqnarray*}\label{J1pre} (\ell)
\lim_w \frac{\rho [2\,\tau \, M(T_w e_0 -e_0)]}{\xi^0_w}=0.
\end{eqnarray*}
and a fortiori
\begin{eqnarray}\label{J1}
(\ell)\lim_w \frac{\rho [2\, \tau \, M(T_w e_0 -e_0)]}{\xi_w}=0,
\end{eqnarray}
taking into account axiom (e) in Definition
\ref{axiomaticconvergence}.
Let $N >0$ be with $|a_r(s)| \leq N$ for any $r=0, \ldots, m$ and
$s \in G$. By (\ref{pi}) and (P1) we get
\begin{eqnarray}\label{cycle1}
P_s(t)=P_s(t)-P_s(s)=\sum_{r=0}^m a_r(s) (e_r(t)-e_r(s))
\end{eqnarray} for each $s$, $t\in G$. By applying 
$T_w$, from (\ref{cycle1}) we obtain
\begin{eqnarray*}
T_w(P_s)(s)=T_w(P_s)(s)-P_s(s)=\sum_{r=0}^m a_r(s) (T_w e_r(s)-e_r(s)),
\end{eqnarray*}
and hence
\begin{eqnarray}\label{cycle3}
|T_w(P_s)(s)|=|T_w(P_s)(s)-P_s(s)|\leq N\sum_{r=0}^m |T_w e_r(s)-e_r(s)|,
\end{eqnarray} for every $s\in G$ and $w\in W$.
By applying the modular $\rho$ and taking into
account $Q$-quasi semiconvexity, from (\ref{cycle3}) we have
\begin{eqnarray}\label{eder}
\rho[ 2 \, \tau \,
{C_1}^{-1} C_2 (T_w P_{(\cdot)}) (\cdot)] \nonumber &\leq&
\sum_{r=0}^m \rho [ 2 \, \tau \,
{C_1}^{-1} C_2 (m+1)  N  (T_w e_r -e_r)] \leq 
\\ &\leq& \sum_{r=0}^m\rho\Bigl[ \frac{\gamma}{Q^2}(T_w e_r - e_r)\Bigr]
 \leq \\ &\leq& Q \, \frac1Q \sum_{r=0}^m\rho \Bigl[ Q \,
\frac{\gamma}{Q}(T_w e_r -e_r)\Bigr]=\sum_{r=0}^m \rho [\gamma (T_w e_r 
-e_r)] \nonumber \end{eqnarray}
for each $w \in W$.
Since $$(\ell)
\lim_w \frac{\rho[\gamma(T_w e_r - e_r)]}{\xi^r_w}=0
\quad \text{  for  every  } r=0,1,\ldots,m,$$ we get
$$(\ell)\lim_w \frac{\rho[\gamma(T_w e_r - e_r)]}{\xi_w}=0
\quad \text{  for  every  } r=0,1,\ldots,m,$$ and hence
\begin{eqnarray}\label{J2}
(\ell)\lim_w\frac{\rho[ 2 \gamma
{C_1}^{-1} C_2 (T_w P_{(\cdot)}) (\cdot)]}
{\xi_w}=0 ,
\end{eqnarray}
by (\ref{eder}) and taking into account axiom (e) of
\ref{axiomaticconvergence}. From (\ref{modular}),
(\ref{J1}) and (\ref{J2}) we obtain
\begin{eqnarray*}
(\ell)\lim_w \frac{\rho[\gamma(T_w f - f)]}{\xi_w}=0,
\end{eqnarray*}
taking into account axioms (a) and (e) in Definition
\ref{axiomaticconvergence}. This ends the proof. $\quad \Box$
\begin{theorem}\label{ratestriple}
Let $Q$, $(T_w)_w$, $\rho$, $\Xi$ be as in 
Theorem \rm \ref{ratesdouble}\em, $(G,d)$ satisfy
condition \rm H*), \em $\xi^0_w$, $\xi^*_w \in \Xi$, set 
$\xi_w:=\max \{\xi^0_w, \xi^*_w\}$, $w\in W$, and
$\psi(s)(t):=d(s,t)$, $s$, $t\in G$. For every
$f \in {\mathcal C}_c(G)$ and $w\in W$
put $\delta^f_w=\|T_w(\psi)\|$, where 
$\|\cdot\|$ is the sup-norm and the supremum is
taken with respect to the support of $f$.
If $\gamma>0$ satisfies the conditions
\begin{description}
\item[{\rm \ref{ratestriple}.1)}] $\rho[\gamma (T_w e_0 - e_0)]
=o(\xi^0_w)$ and
\item[{\rm \ref{ratestriple}.2)}]
$\rho[\gamma \, \omega(f;\delta^f_w)]=o(\xi^*_w)$,
\end{description} then for each $f\in {\mathcal C}_c(G) $ there is
a positive real number $\tau$ with 
$\rho[\tau \, (T_w f - f)]=o(\xi_w)$.

Moreover, a similar result holds when the symbol
$o$ is replaced by $O$.
\end{theorem}
{\bf Proof:}
Let $f\in {\mathcal C}_c(G)$,
$\displaystyle{M=\sup_{t\in G}|f(t)|}$. Observe
that $\omega(f;\delta)\leq 2 \, M$ for each $\delta >0$.
By the properties
of the modulus of continuity, we get
\begin{eqnarray}\label{beginning}
|f(s)-f(t)|\leq \omega(f;d(s,t)) \leq \Bigl (1 +
\frac{d(s,t)}{\delta} \Bigr) \omega(f;\delta)
\end{eqnarray} for each $\delta >0$ and $s$, $t \in G$.
We claim that $\delta^f_w \in \mathbb{R}$ for every $w\in W$.
Indeed the support of $f$ is (totally) bounded, and so
$\displaystyle{
\sup_{s,t\in G} d(s,t) < +\infty}$. By applying $T_w$
we find an $E_w>0$
with $$\displaystyle{\sup_{s,t\in G} T_w(d(s,t)) \leq E_w \,
T_w (e_0)},$$ getting the claim. Let 
$\delta=\delta^f_w$. By applying
$T_w$, keeping fixed $s$ and
letting $t$ vary in $G$, by (\ref{beginning}),
linearity and monotonicity of $T_w$ we get
\begin{eqnarray}\label{middle}
& & |(T_wf)(s)-f(s)| \leq T_w \Bigl(1 +
\frac{d(s,t)}{\delta} \omega(f;\delta)\Bigr) \nonumber \leq
\\ &\leq& \omega(f;\delta) |(T_w e_0) (s) - e_0 (s) | +
\frac{\omega(f;\delta)}{\delta} T_w(\psi)(s)+
\omega(f;\delta)\, e_0(s)+ \\ &+& M |(T_w e_0) (s)-e_0(s)| \leq
4 \, M (|(T_w e_0) (s)-e_0(s)| +
\omega(f;\delta) ) \nonumber
\end{eqnarray} for each $s\in G$. Let now $\gamma> 0$ be 
as in the hypotheses, and pick $\tau >0$ with
 $$ \tau \leq \frac{\gamma}{8 \, M \, Q^2}.$$
  By applying the modular $\rho$, taking into account 
$Q$-quasi semiconvexity, 
from (\ref{middle}) we obtain
\begin{eqnarray}\label{prefinal}
\rho[\tau(T_wf-f)]&\leq& \rho[8 \,\tau \, M (T_w e_0 - e_0)]
+\rho[8 \,\tau \, M \, \omega(f;\delta)] \leq  \nonumber\\ &\leq&
\nonumber \,  \rho\Bigl[\frac{\gamma}{Q^2} (T_w e_0 -
e_0)\Bigr] + \rho[\gamma \,  \omega(f;\delta)]) \leq
\\ &\leq&  Q \, \frac1Q \rho \Bigl[ Q \,
\frac{\gamma}{Q}(T_w e_0 -e_0)\Bigr]+
\rho[\gamma \,  \omega(f;\delta)])=\\&=&\rho [\gamma (T_w e_0 
-e_0)] + \rho[\gamma \,  \omega(f;\delta)]). \nonumber
\end{eqnarray} Since $$(\ell)\lim_w \frac{\rho[\gamma (T_w e_0 -e_0)]}{\xi^0_w}=(\ell)
\lim_w  \frac{\rho[\gamma \, \omega(f;\delta)]}{\xi^*_w}=0,$$ 
then from
(\ref{prefinal}), since
\begin{eqnarray*}
0 &\leq& \frac{\rho[\tau(T_w f-f)]}{\xi_w}
\leq \\ &\leq& 8 \, M \, Q^2 \,\Bigl(
\frac{\rho[\gamma (T_w e_0 -e_0)]}{\xi^0_w} + 
\frac{\rho[\gamma \, \omega(f;\delta)]}{\xi^*_w}\Bigr) \quad \nonumber \text{for  every  }
w\in W,\end{eqnarray*} and by virtue of
the axioms (a), (e) of Definition
\ref{axiomaticconvergence}, it follows that $$
(\ell)\lim_w \frac{\rho
[\tau(T_w f-f)]} {\xi_w}=0,$$ that is the 
first assertion. The proof of the 
last part is analogous. $\quad \Box$
\section{Examples and applications}
We now deal with filter convergence, noting that this kind
of convergence satisfies the axioms of Definitions
\ref{axiomaticconvergence} (see also \cite{BC})
and \ref{axiomaticlimsupliminf}
(see also \cite[Theorems 3 and 4]{DEMIRCI}).
\begin{definitions} \label{prel}
\rm (a) Let $W$ be any abstract infinite set.
A nonempty family ${\mathcal F}$ of subsets of
$W$ is a \textit{filter} of $W$ iff
$\emptyset \not \in {\mathcal F}$,
$A \cap B \in {\mathcal F}$ whenever $A$, $B \in {\mathcal F}$ and for every
$A \in {\mathcal F}$ and $B \supset A$ we get $B \in {\mathcal F}$.

(b) If $W=(W, \succeq)$ is a directed set, then for each $w\in W$, set
$M_w:=\{z\in W: z \succeq w\}$.
A filter ${\mathcal F}$ of $W$ is said to be \textit{free} iff
$M_w\in{\mathcal F}$ for every $w\in W$.

(c) A free filter ${\mathcal F}$
of $W$ is called an \textit{ultrafilter} iff for every
set $A \subset W$ we get that either $A \in {\mathcal F}$ or
$W \setminus A \in {\mathcal F}$.

(d) Let ${\mathcal F}$ be a free filter of $W$.
A net $(x_w)_w$ is said to be \textit{${\mathcal F}$-convergent} to
a real number $x$ iff
\begin{eqnarray}\label{filterconv}
\{w\in W: |x_w-x| \leq \varepsilon \} \in {\mathcal F}
\quad \text{  for every  }  \varepsilon >0,
\end{eqnarray}
and in this case
we write $\displaystyle{({\mathcal F})\lim_w x_w=x}$.

(e) Let $\underline{\bf{x}}=(x_w)_w$ be a net in
$\mathbb{R}$, and set
$$A_{\underline{\bf{x}}}= \{a \in \erre: \{w \in W: x_w \geq a \} \not \in {\mathcal F} \}, $$ $$
 B_{\underline{\bf{x}}}= \{b \in \erre: \{w \in W: x_w \leq b \} \not \in {\mathcal F} \}.$$
The \textit{${\mathcal F}$-limit superior} of $(x_w)_w$ is given by
\begin{eqnarray}\label{lub} ({\mathcal F})\limsup_w \, x_w = \left\{
\begin{array}{ll}
\sup B_{\underline{\bf{x}}}, & \text{if } B_{\underline{\bf{x}}}
\neq \emptyset, \\ -\infty, & \text{if } B_{\underline{\bf{x}}}
 = \emptyset. \end{array} \right.
\end{eqnarray}
The \textit{${\mathcal F}$-limit inferior} of $(x_w)_w$ is
\begin{eqnarray}\label{glb} ({\mathcal F})\liminf_w \, x_w = \left\{
\begin{array}{ll}
\inf A_{\underline{\bf{x}}}, & \text{if } A_{\underline{\bf{x}}}
\neq \emptyset, \\ +\infty, & \text{if } A_{\underline{\bf{x}}}
 = \emptyset \end{array} \right.
\end{eqnarray} (see also \cite{DEMIRCI}).
\end{definitions}
Some examples frequently used in the literature
are $(W, \succeq)$ = $(\mathbb{N}, \geq)$,
$W \subset [a, w_0[ \subset \mathbb{R}$
with the usual order,
where $w_0 \in \mathbb{R} \cup \{ + \infty \}$ is
a limit point of $W$, or $(W, \succeq)=(\mathbb{N}^2,
\geq)$ = $(\mathbb{N} \times \mathbb{N}, \geq$), where
in $\mathbb{N}^2$ the symbol $\geq$ denotes the usual
componentwise order (see also \cite{BMV}).
\vspace{3mm}

(f) The \textit{Fr\'{e}chet filter} is the
filter ${\mathcal F}_{\text{cofin}}$ of all subsets of $\mathbb{N}$ 
whose complement is finite. Observe that the limit,
limit superior
and limit inferior with respect to ${\mathcal F}_{\text{cofin}}$ coincide with the usual
ones (see also \cite{DEMIRCI}).
\begin{remark}\rm 
It is not difficult to check that, if ${\mathcal F}$ is any fixed free filter
of $W$, in the ${\mathcal F}$-convergence setting,
given $(x_w)_w$, $(y_w)_w \in \Xi$,
we get $x_w=o(y_w)$ if and only if 
$\{w \in W$: $x_w \leq \varepsilon \, y_w\}\in{\mathcal F}$ 
for every $\varepsilon >0$ and  
$x_w=O(y_w)$ if and only if there is a positive real number $C$ with
$\{w \in W$: $x_w \leq C \, y_w\}\in{\mathcal F}$. So, our 
Korovkin-type theorems about convergence and rates of approximations 
with respect to the axiomatic convergence in Definitions
\ref{axiomaticconvergence} and \ref{axiomaticlimsupliminf}
contain \cite[Lemma 2.4, Theorem 2.5 and Corollaries 2.6-2.8]{duman}.
\end{remark}
In the filter convergence
context, it is possible also to relax the positivity condition
on the involved linear operators. For instance, let
$I$ be a bounded interval of $\mathbb{R}$, ${\mathcal C}^2(I)$ 
(resp. ${\mathcal C}_b^2(I)$) be the
space of all functions defined on $I$, (resp. bounded and) continuous together with
their first and second derivatives, ${\mathcal C}_+:= \{f \in {\mathcal C}_b^2(I):
f \geq 0 \}$, ${\mathcal C}^2_+:= \{f \in {\mathcal C}_b^2(I):
f^{\prime \prime} \geq 0 \}$.

Let $e_r$, $r=1, \ldots, m$
and $a_r$, $r=0, \ldots, m$ be functions in
${\mathcal C}^2_b(I)$, $P_s(t)$, $s,t \in I$, be as in (\ref{pi}), and
suppose that $P_s(t)$ satisfies
the above conditions (P1), (P2) and
\begin{description}
\item[{\rm (P3)}] there
is a positive real number $C_0$ with $P^{\prime \prime}_s(t) \geq 
C_0$ whenever $s$, $t \in I$, where the second derivative is taken with respect to $t$.
\end{description}
Some examples in which all properties (P1), (P2) and (P3) are satisfied
can be found in \cite{bbdm}.

We now state the following Korovkin-type theorem for not necessarily
positive linear operators, in the setting of modular filter 
convergence, whose proof is analogous to that of 
\cite[Theorem 5.2]{bbdmkorovkin}.
\begin{theorem}\label{korovkin3}
Let ${\mathcal F}$ be any free filter of $W$,
$\rho$ be as in Theorem \ref{korovkin1},
$e_r$, $a_r$, $r=0, \ldots, m$ and $P_s(t)$,
$s$, $t \in I$, satisfy properties \rm (P1), (P2) \em and \rm (P3).
\em Assume that $T_w$, $w\in W$ is a net of
linear operators which fulfil property $(\rho)$-$(*)$, and
that $\{w\in W: 
T_w({\mathcal C}_+ \cap {\mathcal C}^2_+) \subset {\mathcal
C}_+  \} \in {\mathcal F}
$. If $(T_w e_r)_w$ is $(\ell)$-modularly
${\mathcal F}$-convergent to $e_r$,
$r=0, \ldots, m$ in $L^{\rho}(I)$, then $(T_w f)_w$
is $(\ell)$-modularly ${\mathcal F}$-convergent to $f$ in $L^{\rho}(I)$ for
each $f \in {\mathcal C}^2_b(I)$.

If $(T_w e_r)_w$
is $(\ell)$-strongly ${\mathcal F}$-convergent to $e_r$,
$r=0, \ldots, m$ in $L^{\rho}(I)$, then $(T_w f)_w$
is $(\ell)$-strongly 
${\mathcal F}$-convergent to $f$ in $L^{\rho}(I)$ for every $f \in
{\mathcal C}^2_b(I)$.

Furthermore, if $\rho$ is absolutely continuous and $(T_w e_r)_w$
is $(\ell)$-strongly ${\mathcal F}$-convergent to $e_r$,
$r=0, \ldots, m$ in $L^{\rho}(I)$, then $(T_w f)_w$
is $(\ell)$-modularly ${\mathcal F}$-convergent to $f$
in $L^{\rho}(I)$ for every $f \in L^{\rho}(I) \cap {\mathcal D}$
with $f - {\mathcal C}_b(I) \subset X_{{T}}$.
\end{theorem}

\medskip

Other examples of convergences, satisfying 
axioms (a)-(j) in Definitions
\ref{axiomaticconvergence} and \ref{axiomaticlimsupliminf}, are the
\textit{single convergence} and the
\textit{almost convergence} (see also \cite{bdpPL}).

Let $W=\mathbb{N}$. A sequence $(x_n)_n$ is said to 
\textit{singly converge}
(resp. \textit{almost converge}) to $x \in \mathbb{R}$   
iff $$\lim_n \frac{x_{m+1} + x_{m+2} + \ldots + x_{m+n}}{n}=x$$ 
for every $m \geq 0$ (resp. 
uniformly with respect to $m$), where the involved limit 
is the usual one. It is not difficult to 
check that single and almost convergence satisfy
axioms (a)-(j) in Definitions
\ref{axiomaticconvergence} and \ref{axiomaticlimsupliminf}.
In \cite{bdpPL} it is proved that almost convergence is strictly 
stronger than single convergence, and that for every filter
${\mathcal F} \neq {\mathcal F}_{\text{cofin}}$ of $\mathbb{N}$
there are sequences, ${\mathcal F}$-convergent to a point $x_0$ 
but not singly convergent, and a fortiori not almost convergent, to
$x_0$. Moreover an example of a sequence, almost convergent (and
a fortiori singly convergent) to $0$ but not 
${\mathcal F}$-convergent to $0$ for any free filter of $\mathbb{N}$,
is given. Thus, in general, single and almost convergence are not
generated by any free filter.
\vspace{3mm}

We now consider a kind of ``triangular statistical convergence''
investigated in \cite{bbdmo2, dobbm}. 
Let $A=(a_{i,j})_{i,j}$ be a non-negative
two-dimensional infinite matrix and $\Psi:\mathbb{N} \times
\mathbb{N} \to \mathbb{R}$ be a fixed function.
We say that $A$ is a \textit{summability matrix} iff it satisfies 
the following conditions:



\medskip

$(A1)$  $\displaystyle{
\sum_{j\in\mathbb{N},\Psi(i,j)\geq 0}
a_{i,j}\leq 1}$ for each $i\in\mathbb{N}$,

$(A2)$ $\displaystyle{\lim_i 
\sum_{j\in\mathbb{N},\Psi(i,j)\geq 0}
a_{i,j}>0}$,

$(A3)$ $\displaystyle{\lim_i a_{i,j}=0}$ for every $j\in \mathbb{N}$

(see also \cite{KL}).

\bigskip

For every $K\subset \mathbb{N}^2$, set 
$K_i:=\{j\in\mathbb{N}$: $(i,j)\in K$,
$\Psi(i,j)\geq 0 \}$.
The \textit{$\Psi$-$A$-density} of $K$ is given by
\begin{eqnarray}\label{psiA}
\delta_A^{\Psi}(K):=\lim_i\sum_{j\in K_{i}}a_{i,j},
\end{eqnarray}
provided that the limit on the right hand exists in $\mathbb{R}$.

It is not difficult to see that the $\Psi$-$A$-density satisfies
the following properties:

\medskip
$(D1)$ $\delta_A^{\Psi}(\mathbb{N}^{2})>0$.

$(D2)$ If $K\subset H$, then $\delta_{A}^{\Psi}(K)\leq 
\delta_A^{\Psi}(H)$.


$(D3)$ If $\delta_A^{\Psi}(K)=\delta_A^{\Psi}(H)=0$, then 
$\delta_A^{\Psi}(K\cup H)=0$

(see also \cite{bbdmo2}). Observe that from $(D1)$-$(D3)$ it 
follows that the family
\begin{eqnarray}\label{fapsi}
{\mathcal F}_A^{\Psi}:=\{ K \subset \mathbb{N}^2:
\delta_A^{\Psi}(\mathbb{N}^2 \setminus K)=0\} 
\end{eqnarray}
is a filter of $\mathbb{N}^2$.
In order to show that ${\mathcal F}_A^{\Psi}$ is free, it will
be enough to check that, if $p$, $q\in\mathbb{N}$ and 
$K:=\{(i,j)\in \mathbb{N}^2$: $i \geq p$, $j \geq q\}$, then 
$\delta_A^{\Psi}(\mathbb{N}^2 \setminus K)=0$. Indeed, we get
$\mathbb{N}^2 \setminus K \subset K^{(1)} \cup K^{(2)}$, where
$$\displaystyle{K^{(1)}=\bigcup_{j=1}^{q-1}(\mathbb{N} \times 
\{j\})}, \quad \displaystyle{K^{(2)}=\bigcup_{i=1}^{p-1}
(\{i\} \times \mathbb{N})}. $$ From $(A3)$ it follows that 
$$\delta_A^{\Psi}(K^{(1)})=\lim_i \sum_{j\in[1,q-1],
\Psi(i,j)=0} a_{i,j}\leq \lim_i \sum_{j=1}^{q-1} a_{i,j}=0.$$
Furthermore, observe that $\delta_A^{\Psi}(K^{(2)})=0$, since 
$K^{(2)}_i=\emptyset$ for every $i \geq p$. Hence, 
$\delta_A^{\Psi}(\mathbb{N}^2 \setminus K)=0$, that is the claim.
\begin{definition}
\label{def1} \rm Let $A=(a_{i,j})_{i,j}$ be a
summability matrix.
The double sequence $(x_w)_w$ is said to 
\textit{$\Psi$-$A$-statistically converge}
to a real number $x$ iff $\displaystyle{
({\mathcal F}_A^{\Psi})\lim_w x_w=x}$, that is iff
for every $\varepsilon >0$ we get
\begin{eqnarray*}
\lim_i \sum_{j\in K_i(\varepsilon)}a_{i,j}=0,
\end{eqnarray*}
where $K_i(\varepsilon) =\{j\in \mathbb{N}$:
$\Psi(i,j)\geq 0$, $|x_w-x| \geq \varepsilon \}$, and we write
$st_{A}^{\Psi}$-$\displaystyle{\lim_i x_w=x}$.

By $st_{A}^{\Psi}$-$\displaystyle{\limsup_i x_w}$ and
$st_{A}^{\Psi}$-$\displaystyle{\liminf_i x_w}$ we denote the
quantities $\displaystyle{({\mathcal F}_A^{\Psi})\limsup_w 
x_w}$ and $\displaystyle{({\mathcal F}_A^{\Psi})\liminf_w 
x_w}$, defined as in (\ref{lub}) and (\ref{glb}), respectively.
\end{definition} 

Let $(a_{i,j})_{i,j}$ be defined by
\begin{eqnarray*}\label{cesaro2}
a_{i,j}:= \left\{\ \begin{array}{ll} \frac{1}{i^2} & \text{  if
} j \leq i^2, \\ \\0 & \text{  otherwise,} \end{array}\right.
\end{eqnarray*} put $\Psi (i,j) = i-j$, $i$, $j\in\mathbb{N}$,
and pick any double sequence $(x_{i,j})_{i,j}$ in $\mathbb{R}$.
For every $\varepsilon> 0$ we get  
$K_i(\varepsilon)$ := $\{j\in \mathbb{N}$:
$j\leq i$, $x_{i,j}\geq \varepsilon \}\subset\{j\in \mathbb{N}$:
$j\leq i\}$. Thus we obtain
\begin{eqnarray}
\lim_i \sum_{j\in K_i(\varepsilon)}a_{i,j} \leq
\lim_i \sum_{j\leq i} \frac{1}{i^2}= 
\lim_i \frac{1}{i}=0,
\end{eqnarray}
and thus $(x_{i,j})_{i,j}$ 
$\Psi$-$A$-statistically converges to $0$. We get
$\displaystyle{\lim_i
\sum_{j\in\mathbb{N},\Psi(i,j)\geq 0}
a_{i,j}=0}$, and so condition $(A2)$ is not fulfilled. 
Note that, in this case, the class ${\mathcal F}_A^{\Psi}$ defined as 
in (\ref{fapsi}) is not a filter, because it coincides
with the family of all subsets of $\mathbb{N}^2$.

\begin{remarks} \rm (a) If we take $\Psi (i,j) = i-j$, $i$, 
$j\in\mathbb{N}$, then we obtain the notion of \textit{triangular 
$A$-statistical convergence}.

(b) If $A=C_1$ is the Ces\`{a}ro matrix, defined by setting
\begin{eqnarray*}\label{cesaro}
a_{i,j}:= \left\{\ \begin{array}{ll} \frac1i & \text{  if
} j \leq i, \\ \\0 & \text{  otherwise,} \end{array}\right.
\end{eqnarray*}
then the $\Psi$-$A$-statistical convergence can be viewed as an 
extension of the classical statistical convergence
(see also \cite{bbdmo2}).
\end{remarks}
We now claim that, in general, filter convergence in 
$\mathbb{N}^2$ is not equal to $\Psi$-$A$-statistical
convergence, that is there exists some filter ${\mathcal F}$
of $\mathbb{N}^2$ such that, for every summability matrix $A$,
there is a set $K \in {\mathcal F} \setminus {\mathcal F}_A^{\Psi}$.
Fix arbitrarily a summability matrix $A=(a_{i,j})_{i,j}$. By 
$(A1)$ and $(A2)$
there exist a real number $B_0 \in (0,1]$ and an infinite subset 
$S\subset \mathbb{N}$ with 
$\displaystyle{\sum_{j\in\mathbb{N},\Psi(i,j)\geq 0}
a_{i,j} \geq B_0}$ whenever $i\in S$.  Let $S_1$ and $S_2$ be
two disjoint infinite subsets of $S$, with $S=S_1\cup S_2$. 
Let $K:=\{(i,1)$: $i\in\mathbb{N}\}$ $\cup$ $\{(i,j): i \in S_1$,
$\Psi(i,j)\geq 0\}$. Taking into account $(A3)$, we get 
\begin{eqnarray}\label{zero} 0 &\leq& \nonumber
\displaystyle{\liminf_{i\in \mathbb{N} \setminus S_1} 
\sum_{j\in\mathbb{N}: (i,j)\in K, \Psi(i,j)\geq 0} a_{i,j} \leq
\limsup_{i\in \mathbb{N} \setminus S_1} 
\sum_{j\in\mathbb{N}: (i,j)\in K, \Psi(i,j)\geq 0} a_{i,j}}
\leq \\ &\leq& \displaystyle{\limsup_i a_{i,1}=\lim_i a_{i,1}=0}.
\end{eqnarray} 
Thus all inequalities in  (\ref{zero}) are equalities, and in particular
we have \begin{eqnarray}\label{zero1}
\lim_{i\in \mathbb{N} \setminus S_1} 
\sum_{j\in\mathbb{N}: (i,j)\in K, \Psi(i,j)\geq 0} a_{i,j}=0.
\end{eqnarray}
Since $\displaystyle{\limsup_{i\in S_1}
\sum_{j\in\mathbb{N}: (i,j)\in K,
\Psi(i,j)\geq 0} a_{i,j} \geq B_0}$, from this and (\ref{zero1}) it 
follows that $\mathbb{N}^2 \setminus 
K \not \in {\mathcal F}_A^{\Psi}$. From this and $(A2)$ it follows
that $K \not \in {\mathcal F}_A^{\Psi}$. Thus we get that any
ultrafilter ${\mathcal F}$ of $\mathbb{N}$
contains at least a set not belonging to 
${\mathcal F}_A^{\Psi}$, because it contains either $K$ or
$\mathbb{N}^2 \setminus K$. This proves the claim. So, our results
extend \cite[Theorems 1-4]{bbdmo2} and 
\cite[Theorems 3-5]{dobbm}.

\bigskip

We now give some applications to Mellin-type convolution operators
and Kan\-to\-ro\-vich-type discrete operators. 
\begin{examples} \rm (a)
We deal with a direct extension to the multivariate case of the
classical one-dimensional moment kernel (see also 
\cite{bm2008, bcs}).

Let $W=[1,+\infty[$ or $W=\mathbb{N}$, $G=[0,1]^N$
equipped with the usual topology, $\rho$ be
as in Theorem \ref{ratesdouble}.
Let ${\mathcal F}$ be a free filter of $W$, containing a subset $F
\subset W$ such that $W \setminus F$ is infinite. An example is 
the filter ${\mathcal F}_{\text{st}}$ of all subsets of $\mathbb{N}$ 
having asymptotic density one, and $\mathbb{N} \setminus F$ 
is the set of all perfect squares or the set of all prime numbers 
(see also \cite{bdrates}).
For every $w\in W$ and $\textbf{t}$ =$(t_1$, $t_2, \ldots, t_N)\in G$,
let $K_w(\textbf{t})=(w+1)^N t_1^w \cdot \ldots \cdot t_N^w$ if 
$w \in F$ and $K_w(\textbf{t})=(w+1)^{N+1} t_1^w \cdot \ldots \cdot t_N^w$ if $w \in W \setminus F$. For
$f\in {\mathcal C}(G)$ and $\textbf{s}$ = $(s_1$, $s_2, \ldots, s_N)\in G$ set
$$(M_w f)(\textbf{s})=\int_{G} K_w(\textbf{t}) f(\textbf{s} 
\textbf{t}) \, d\textbf{t},$$ where $\textbf{s}\textbf{t}$ = $(s_1 t_1, 
s_2 t_2, \ldots, s_N t_N)$ and $d\textbf{t}=dt_1 \, dt_2\, \ldots \, 
dt_N$. For each $\textbf{t}\in G$, set $e_0(\textbf{t})=1$, 
$e_r(\textbf{t})=t_r$, $r=1, \ldots, N$, and $\displaystyle{
e_{N+1}(\textbf{t})=\sum_{r=1}^N t_r^2}$. It is easy to see that 
(P1), (P2) and (P3) are satisfied.
We get
\begin{eqnarray*} 
\int_G K_w(\textbf{t}) d\textbf{t}=(w+1)^N  \left(\int_0^1 t_1^w \, 
dt_1 \right)
\cdot \ldots \cdot \left(\int_0^1 t_N^w \, dt_N \right)=1
\end{eqnarray*}
if $w\in F$, and \begin{eqnarray*} 
\int_G K_w(\textbf{t}) d\textbf{t}=(w+1)^{N+1} \left(\int_0^1 
t_1^w \, dt_1 \right)
\cdot \ldots \cdot \left(\int_0^1 t_N^w \, dt_N \right)=w+1
\end{eqnarray*}
if $w \in W \setminus F$.
Hence for every $\textbf{s}\in G$ we get $(M_w e_0)
(\textbf{s})=e_0(\textbf{s})=1$ if $w\in F$ and $(M_w e_0)
(\textbf{s})-e_0(\textbf{s})=w$ if $w\in W \setminus F$.
From this it follows that the operators $M_w$, $w\in W$, 
does not satisfy the classical Korovkin theorem. 
However we have
\begin{eqnarray*}
& & |(M_w e_r)(\textbf{s})-e_r(\textbf{s})|
\leq \int_{G} K_w(\textbf{t})(1-t_r) \, d\textbf{t}=\\ &=&
(w+1) \Bigl( \int_0^1 t_1^w(1-t_1) \, dt_1 \Bigr) \cdot
\ldots \cdot \Bigl( \int_0^1 t_N^w(1-t_N) \, dt_N \Bigr)=\\&=& (w+1)
\Bigl( \int_0^1 t_r^w(1-t_N) \, dt_r \Bigr)
= 1-\frac{w+1}{w+2}=\frac{1}{w+2}=O\Bigl(\frac1w\Bigr)
\end{eqnarray*} for every $r=1$, $2,\ldots, N$,
$w\in F$ and $\textbf{s}\in G$.  
Analogously it is possible to see that
\begin{eqnarray*}
& & |(M_w e^2_r)(\textbf{s})-e^2_r(\textbf{s})|
\leq \frac{2}{w+3}=O\Bigl(\frac1w\Bigr)
\end{eqnarray*} whenever $r=1$, $2,\ldots, N$, 
$w\in F$ and $\textbf{s}\in G$. From this it 
follows that $$|(M_w e_{N+1})(\textbf{s})-e_{N+1}(\textbf{s})|=
O\Bigl(\frac1w\Bigr)$$ for $w\in W$ and $\textbf{s}\in G$.
Thus there exists a positive real number $C_*$ with
\begin{eqnarray}\label{premodularestimate}
|(M_w e^2_r)(\textbf{s})-e^2_r(\textbf{s})|
\leq \frac{C_*}{w}
\end{eqnarray} for any $r=0,1$, $2,\ldots, N+1$, 
$w\in F$ and $\textbf{s}\in G$. Let $Q\geq 1$ be a constant
related with $Q$-quasi semiconvexity of $\rho$. 
By applying the modular $\rho$, we get
\begin{eqnarray}\label{Qsemi}
\rho\Bigl[\frac{1}{Q^2}(M_w e_r - e_r)\Bigr] \leq \rho 
\Bigl[\frac{C_*}{w \, Q^2} \Bigr] \leq \frac{C_*}{w} \, \rho[1] 
\end{eqnarray} for each $r=0,1$, $2,\ldots, N+1$ and  
$w\in F$. Thus, thanks to strong finiteness of $\rho$,
all the hypotheses of Theorem \ref{ratesdouble}
are satisfied. So for every $f\in 
{\mathcal C}_c(G)\cap$ Lip$(G)$ there is 
$\tau > 0$ with $\displaystyle{
\rho[\tau(M_w f - f)]=O\Bigl(\frac1w\Bigr)}$ . 

\medskip

(b) We consider bivariate Kantorovich-type operators. Let 
$W=\mathbb{N}$, ${\mathcal 
F}\neq{\mathcal F}_{\rm{cofin}}$
and $H$ be an infinite set with $\mathbb{N} \setminus H \in 
{\mathcal F}$. Note that $H$ does exist, since 
${\mathcal F}\neq{\mathcal F}_{\rm{cofin}}$. Let $W=\mathbb{N}$,
$G=[0,1]^2$ be endowed with the usual topology, $\widetilde{\Phi}$ 
be as in Example \ref{orlicz}, $\varphi \in \widetilde{\Phi}$ be with 
$\displaystyle{\liminf_{u \to\infty} \frac{\varphi(u)}{u} >0}$,
$\rho=\rho^{\varphi}$ be as in (\ref{modularorlicz}).
Proceeding as in \cite{belen}
(see also \cite{bbdmkorovkin}),
for every locally integrable function
$f \in L^0(G)$, $n \in \mathbb{N}$ and $x$, $y \in [0,1]$ set
$$P_n(f)(x,y)=(n+1)^2 \sum_{k,j=0,\ldots,n, k+j \leq n} p_{n,k,j} (x,y)
\int_{k/(n+1)}^{(k+1)/(n+1)}\int_{j/(n+1)}^{(j+1)/(n+1)} f(u,v) \, du \, dv,$$ where
$$p_{n,k,j}=\frac{n!}{k!j!(n-k-j)!} x^k y^j (1-x-y)^{n-k-j}, \quad k,j \geq 0, x, y\geq 0,
x+y \leq 1.$$
Let $(s_n)_n$ be the sequence defined by setting
\begin{eqnarray*}
s_n=\left\{ \begin{array}{ll} 1 & \text{if  } n \in \enne \setminus H, \\
0& \text{if  } n \in H. \end{array} \right.
\end{eqnarray*}
For every $n \in \mathbb{N}$ and $x$, $y \geq 0$ with $x+y \leq 1$, 
set $$P^*_n(f) (x,y)= s_n \, P_n(f)(x,y).$$
For each $t_1$, $t_2 \in [0,1]$, set $e_0(t_1,t_2)=1$, 
$e_1(t_1,t_2)=t_1$, 
$e_2(t_1,t_2)=t_2$, $e_3(t_1,t_2)=t_1^2+t_2^2$ (see also 
\cite[Example 4.5(b)]{bbdmkorovkin}).
By proceeding analogously as in \cite{belen}, it is possible to find a
$\gamma > 0$ with $\displaystyle{
\rho[\gamma \, (P^*_n e_i - e_i)] =O\Bigl(\frac{1}{n}\Bigr)}$.
By Theorem \ref{ratesdouble}, for every  
$f\in {\mathcal C}_c(G) \cap \text{Lip}(G)$ there is a positive real
number $\tau$ with 
$\displaystyle{\rho[\tau \, (P^*_n f - f)]=O
\Bigl(\frac{1}{n}\Bigr)}$. 
\end{examples}
{\bf Acknowledgements:} Supported by Universities of Perugia 
and Athens.

\end{document}